\documentclass[10pt]{amsart}
\usepackage{amsfonts,amssymb}
\def\qq{\mathbb{Q}}
\def\rr{\mathbb{R}}
\def\cc{\mathbb{C}}
\def\zz{\mathbb{Z}}

\def\hh{\mathcal{H}}
\def\xx{\mathcal{X}}

\def\ww{\mathcal{W}}
\def\bigone{\mathbf{1}}
\def\Ardeg{\widehat{\deg}}
\def\dual{{\omega}_{\xx / B}}
\def\height{\Ardeg \det p_* \dual}

\def\difftials{H^{0}(X,\Omega_{X}^{1})}
\def\imtau{\mathrm{Im}\tau}
\def\mm{\mathcal{M}}

\def\isom{{\buildrel \sim \over \longrightarrow}}

\theoremstyle{plain}
\newtheorem{theorem}{Theorem}[section]
\newtheorem{lemma}[theorem]{Lemma}
\newtheorem{proposition}[theorem]{Proposition}
\newtheorem{corollary}[theorem]{Corollary}

\theoremstyle{remark}

\theoremstyle{definition}

\theoremstyle{definition}
\newtheorem{definition/proposition}[theorem]{Definition/Proposition}

\numberwithin{equation}{section}

\begin{document}

\title{Arakelov invariants of Riemann surfaces}
\author{Robin de Jong}
\address{University of Amsterdam, The Netherlands}
\email{rdejong@science.uva.nl}

\begin{abstract}
We derive explicit formulas for the Arakelov-Green function and
the Faltings delta-invariant of a Riemann surface. A numerical
example illustrates how these formulas can be used to calculate
Arakelov invariants of curves.
\end{abstract}

\date{\today}

\maketitle

\section{Introduction}

The Arakelov-Green function and Faltings' delta-invariant are
fundamental invariants attached to Riemann surfaces \cite{ar},
\cite{fa}. However, they are defined in a quite implicit way. It is
therefore natural to ask for explicit formulas for these invariants,
and indeed in many cases such explicit formulas are known. For example,
in \cite{fa} the case of elliptic curves is treated in detail,
and in \cite{bo}, \cite{bmmb} we find explicit results dealing with the case of
Riemann surfaces of genus 2. In higher genera there only seem to be 
some scattered
results; for example, \cite{gu3} treats a certain family of
plane quartic curves, and in \cite{abbull}, \cite{michull}
the modular curves $X_0(N)$ are studied from the point of view of
Arakelov theory. The purpose of the
present note is to give general explicit formulas for the
Arakelov-Green function and the delta-invariant that make it possible
to calculate these invariants efficiently. We have included an explicit
numerical example at the end of this note to illustrate the use of our
formulas for computations.

We now describe our results. Let $X$ be a
compact and connected Riemann surface of genus $g>0$. We recall
from \cite{ar} and \cite{fa} that $X$ carries a canonical
(1,1)-form $\mu$, giving rise to a Green-function $G : X \times X
\to \rr_{\geq 0}$ and a canonical structure of metrised line
bundle on the holomorphic cotangent bundle $\Omega_X^1$ and the
line bundles $O_X(D)$ associated to a divisor $D$ on $X$. The line
bundle $O(\Theta)$ on $\mathrm{Pic}_{g-1}(X)$ admits a metric $\|
\cdot \|_{\Theta}$ with $\|s\|_{\Theta}=\|\vartheta\|$, where $s$
is the canonical section of $O(\Theta)$ and where $\|\vartheta\|$
is the function defined on \cite{fa}, p.~401.

Our first result deals with the Arakelov-Green function $G$. It
has been observed by some authors (see the remarks on \cite{jo},
p.~229) that for a generic point $P \in X$ there exists a
constant $c=c(P)$ depending only on $P$ such that for all $Q \in
X$ we have $G(P,Q)^g = c(P) \cdot \|\vartheta\|(gP-Q)$. Our
contribution is that we make the dependence on $P$ of the constant
$c(P)$ clear. Our formula involves the divisor $\ww$ of
Weierstrass points on $X$. Recall that this is a divisor of degree
$g^3-g$ on $X$, given as the divisor of a Wronskian differential
formed out of a basis of the holomorphic differentials
$\difftials$. For each point $P \in X$, the multiplicity of $P$ in
$\ww$ is given by a weight $w(P)$, that can also be calculated by
means of the classical gap sequence at $P$.

Let $S(X)$ be the invariant defined by the formula \[ \log S(X) :=
- \int_X \log \|\vartheta\|(gP-Q) \cdot \mu(P) \] for any $Q \in
X$. We will see
later that the integrand has logarithmic
singularities only at the Weierstrass points of $X$, which are
integrable. Hence the integral is well-defined. Also we will see later
that the integral does not depend on the choice of $Q$.
As an example, consider the case $g=1$ and write
$X=\cc/\zz+\tau \zz$, with $\tau$ in the complex upper half plane.
A calculation (see for example \cite{la2}, p. 45) shows that in this
case
\[ \log S(X) = - \log ( (\imtau)^{1/4} |\eta(\tau)| ) \, , \]
where $\eta(\tau)$ is the usual Dedekind eta-function given by
$\eta(q) := q^{1/24} \prod_{n=1}^\infty (1-q^n)$, with
$q=\exp(2\pi i \tau)$.

The invariant $S(X)$ appears as a normalisation constant in the
formula that we propose for the Arakelov-Green function.
\begin{theorem} \label{Green} Let $P,Q \in X$ with $P$ not a
Weierstrass point. Then the formula \[ G(P,Q)^g = S(X)^{1/g^2}
\cdot \frac{ \| \vartheta \|(gP-Q) }{ \prod_{W \in \ww} \|
\vartheta \|(gP-W)^{1/g^3}} \] holds. Here the Weierstrass points
are counted with their weights.
\end{theorem}
For $P$ a Weierstrass point, and $Q \neq P$,
both numerator and denominator in the
formula of Theorem \ref{Green} vanish with order $w(P)$, the
weight of $P$. The formula remains true also in this case,
provided that we take the leading coefficients of the appropriate
power series expansions about $P$ in both numerator and
denominator. Note that apart from the normalisation term involving
$S(X)$, the Arakelov-Green function can be expressed in terms of
certain values of the $\|\vartheta\|$-function. These values are
very easy to calculate numerically. The (real)
2-dimensional integral involved in computing $S(X)$ is harder to
carry out in general, but it is still not difficult.

Other ways of expressing the Arakelov-Green function in terms of
quantities associated to $X$ and $\mu$ have been given, for
instance one might use the eigenvalues and eigenfunctions of the
Laplacian (see \cite{fa}, Section 3), or one might use abelian differentials of
the second and third kind (see \cite{la2}, Chapter II). There is
also a closed formula due to Bost \cite{bo}
\[ \log G(P,Q) = \frac{1}{g!} \int_{\Theta+P-Q} \log \|\vartheta\|
\cdot \nu^{g-1} + A(X) \, , \] expressing the Arakelov-Green
function in terms of an integral over the translated theta
divisor. Here $\nu$ is the canonical translation-invariant
(1,1)-form on $\mathrm{Pic}_{g-1}(X)$, and the quantity $A(X)$ is
a certain normalisation constant, perhaps comparable to our
$S(X)$.

One of our motives for finding a new explicit formula was the need
to have a formula that makes the efficient calculation of the
Arakelov-Green function possible. The other approaches that we
mentioned are perhaps less suitable for this objective. For
instance, the formula given by Bost involves a (real)
$2g-2$-dimensional integral over a region which seems not easy to
parametrise. Also, for each new pair of points $(P,Q)$ one has to
calculate such an integral again, whereas in our approach one only
has to calculate a certain integral once.

Our second result deals with Faltings' delta-invariant
$\delta(X)$. It is the constant appearing in the following
theorem, due to Faltings (\emph{cf.} \cite{fa}, p.~402). \begin{theorem}
\label{Faltings} (Faltings) There is a constant $\delta=\delta(X)$
depending only on $X$ such that the following holds. Let
$\{\omega_1,\ldots,\omega_g\}$ be an orthonormal basis of
$\difftials$ provided with the hermitian inner product
$(\omega,\eta) \mapsto \frac{i}{2} \int_X \omega \wedge
\overline{\eta}$. Let $P_1,\ldots,P_g,Q$ be generic points on $X$.
Then the formula \[ \| \vartheta \| (P_1+ \cdots + P_g - Q) = \exp
(-\delta(X)/8) \cdot \frac{ \| \det \omega_k(P_l)
\|_{\mathrm{Ar}}}{ \prod_{k<l} G(P_k,P_l)} \cdot
 \prod_{k=1}^{g} G(P_k,Q)  \] holds.
\end{theorem}
The significance of the delta-invariant is that it appears as an
archimedean contribution in the so-called Noether formula
\cite{fa}, \cite{mb} for arithmetic surfaces. When viewed as a
function on the moduli space $\mm_g$ of curves of genus $g$, the
value $\delta(X)$ can be seen as the minus logarithm of the
distance of the class of $X$ to the Deligne-Mumford boundary of
$\mm_g$. This interpretation is supported by the Noether formula.

Let $\Phi : X \times X \to \mathrm{Pic}_{g-1}(X)$ be the map
sending $(P,Q)$ to the class of $(gP-Q)$. For a fixed $Q \in X$,
let $i_Q : X \to X \times X$ be the map sending $P$ to $(P,Q)$,
and put $\phi_Q := \Phi \cdot i_Q$. Define the (fractional) line
bundle $L_X$ by \[ \begin{aligned} L_X := \left( \bigotimes_{W \in
\ww} \phi_W^*  \left(  O(\Theta) \right) \right)&^{\otimes
(g-1)/g^3}  \otimes_{O_X} \left( \Phi^*(O(\Theta))|_{\Delta_X}
\otimes_{O_X} \Omega_X^{\otimes g} \right)^{\otimes -(g+1)}
\otimes_{O_X}
\\ &\otimes \left( \Omega_X^{\otimes g(g+1)/2} \otimes_{O_X} \left(
\wedge^g \difftials \otimes_\cc O_X) \right)^\lor \right)^{\otimes
2} \, .
\end{aligned}
\] We have then the following theorem.
\begin{theorem} \label{delta}
The line bundle $L_X$ is canonically trivial. Let
$T(X)$ be the norm of the canonical trivialising section of $L_X$.
Then the formula \[  \exp(\delta(X)/4) = S(X)^{-(g-1)/g^2} \cdot
T(X) \] holds.
\end{theorem}
Despite appearances to the contrary, the invariant $T(X)$ admits a
very concrete description. In Proposition \ref{Twithfirstorder}
below we will see that the computation of $T(X)$ only involves
elementary operations on special values of $\|\vartheta\|$ and of
the $\|J\|$-function, a function introduced by Gu\`ardia in
\cite{gu}. These special values are easy to calculate numerically.
The significance of Theorem \ref{delta} is then that we have
reduced the calculation of $\delta(X)$ to the calculation of two
new invariants $S(X)$ and $T(X)$, the former involving a (real)
2-dimensional integral over the surface $X$, the latter being
elementary to calculate.

It seems an important problem to relate the invariants $S(X)$ and
$T(X)$ to more classical invariants. In Theorem
\ref{hyperelliptic} we state a result that does this for $T(X)$
with $X$ a hyperelliptic Riemann surface.

The plan of this note is as follows. First in Section \ref{proofs}
we give the proofs of Theorems \ref{Green} and \ref{delta}. The major
idea will be to give Arakelov-theoretical versions of
classical results on the divisor of Weierstrass points. In Section
\ref{application} we will give some applications of our results in
the Arakelov intersection theory of arithmetic surfaces. We derive
a lower bound for the self-intersection of the relative dualising
sheaf, and we give a formula for the self-intersection of a point.
In Section \ref{example} we give a numerical example in the spirit
of \cite{bmmb}, calculating the Arakelov invariants of an
arithmetic surface associated to a certain hyperelliptic curve of
genus 3 and defined over $\qq$.

Our inspiration to study Weierstrass points in order to obtain
results in Arakelov theory stems from the papers \cite{ar2},
\cite{bu} and \cite{jo}. Especially the latter paper has been
useful. For example, our formula for the delta-invariant in Theorem \ref{delta}
is closely related to the formula from Theorem 2.6 of 
that paper. Our improvement on that formula is perhaps that we give 
an explicit splitting
of the delta-invariant in a new invariant $S(X)$ involving an
integral, and a new invariant $T(X)$ which is purely `classical'.
These invariants seem to be of interest in their own 
right, \emph{cf.} also our remarks at the end of Section \ref{proofs}.

\section{Proofs} \label{proofs}

We start by recalling the definitions of the (1,1)-form $\mu$, the
Arakelov-Green function $G$ and the canonical metric on
$\Omega_X^1$. The (1,1)-form $\mu$ is given by $\mu = \frac{i}{2g}
\sum_{k=1}^g \omega_k \wedge \overline{\omega}_k$, where
$\{\omega_1,\ldots,\omega_g\}$ is an orthonormal basis of the
holomorphic differentials $\difftials$ provided with the hermitian
inner product $(\omega,\eta) \mapsto \frac{i}{2} \int_X \omega
\wedge \overline{\eta}$.

The Arakelov-Green function $G$ is
the unique function $ X \times X \to \rr_{\geq 0}$ such that
the following three properties hold:
\begin{itemize} \item[(i)] $G(P,Q)^2$ is $C^\infty$ on $X \times X$
and $G(P,Q)$ vanishes only at the diagonal $\Delta_X$, with multiplicity
1;
\item[(ii)] for all $P \in X$ we have $\partial_Q
\overline{\partial}_Q
\log G(P,Q)^2 = 2\pi i \mu(Q)$ for $Q \neq P$;
\item[(iii)] for all $P \in X$ we have $\int_X \log G(P,Q) \mu(Q)=0$.
\end{itemize} These
properties imply, by an application of Stokes' theorem, the
symmetry $G(P,Q)=G(Q,P)$ of the function $G$.

The canonical metric $\| \cdot \|_{\mathrm{Ar}}$ on the cotangent
bundle $\Omega_X^1$ is the unique metric that makes the canonical
adjunction isomorphism $O_{X \times X}(-\Delta_X)|_{\Delta_X}
\isom \Omega_X^1$ an isometry, the line bundle $O_{X\times
X}(\Delta_X)$ being given the hermitian metric defined by
$\|\bigone_{\Delta_X}\|(P,Q):=G(P,Q)$.

Next let us recall the Wronskian differential that defines the
divisor of Weierstrass points on $X$. For proofs and more details
we refer to \cite{gun}, pp. 120--128.  Let $\{
\psi_1,\ldots,\psi_g \}$ be a basis of $\difftials$. Let $P$ be a
point on $X$ and let $z$ be a local coordinate about $P$. Write
$\psi_k = f_k \cdot dz$ for $k=1,\ldots,g$. The Wronskian
determinant about $P$ is then the holomorphic function \[
W_z(\psi) := \det \left( \frac{1}{(l-1)!} \frac{ d^{l-1} f_k
}{dz^{l-1}} \right)_{1\leq k,l \leq g} \, . \] Let $\tilde{\psi}$
be the $g(g+1)/2$-fold holomorphic differential
\[ \tilde{\psi} := W_z(\psi) \cdot (dz)^{\otimes g(g+1)/2} \, . \]
Then $\tilde{\psi}$ is independent of the choice of the local
coordinate $z$, and extends to a non-zero global section of
$\Omega_X^{g(g+1)/2}$. A change of basis changes the
Wronskian differential by a non-zero scalar factor, so that the
divisor of a Wronskian differential $\tilde{\psi}$ on $X$ is unique: we denote
this divisor by $\ww$, the divisor of Weierstrass points.

The Wronskian differential leads to a canonical sheaf morphism
\[ \left( \wedge^g \difftials \otimes_\cc O_X \right) \longrightarrow
\Omega_X^{g(g+1)/2} \] given by
\[   \xi_1 \wedge \ldots \wedge \xi_g  \mapsto
\frac{\xi_1 \wedge \ldots \wedge \xi_g}{\psi_1 \wedge \ldots
\wedge \psi_g } \cdot \tilde{\psi} \, . \]
This gives a canonical section in $ \Omega_X^{\otimes g(g+1)/2}
\otimes_{O_X} \left( \wedge^g \difftials
\otimes_\cc O_X) \right)^\lor $ whose divisor is $\ww$.
\begin{proposition} \label{normWeier} The canonical isomorphism
\[ \Omega_X^{\otimes g(g+1)/2} \otimes_{O_X} \left( \wedge^g \difftials
\otimes_\cc O_X) \right)^\lor \isom O_X(\ww)  \] has a constant norm
on $X$.
\end{proposition}
\begin{proof} This follows since both sides have the same
curvature form, and the divisors of the canonical sections are equal.
\end{proof}
We shall denote by $R(X)$ the norm of the isomorphism from
Proposition \ref{normWeier}. In more concrete terms we have
$\prod_{W \in \ww} G(P,W) = R(X) \cdot \| \tilde{\omega}
\|_{\mathrm{Ar}}(P)$ for any $P \in X$, where
$\{\omega_1,\ldots,\omega_g\}$ is an orthonormal basis of
$\difftials$, and where the norm of $\tilde{\omega}$ is taken in
the line bundle $\Omega_X^{\otimes g(g+1)/2}$ with its canonical
metric induced from the canonical metric on $\Omega_X^1$. Taking
logarithms and integrating against $\mu(P)$ gives, by property
(iii) of the Arakelov-Green function, the formula $\log R(X) = -
\int_X \log \|\tilde{\omega}\|_{\mathrm{Ar}}(P) \cdot \mu(P)$.

Recall from the Introduction the map $\Phi : X \times X \to
\mathrm{Pic}_{g-1}(X)$ sending $(P,Q)$ to the class of $(gP-Q)$.
A classical result on the divisor of Weierstrass points is that
the equality of divisors
\[ \Phi^*(\Theta) = \ww \times X + g \cdot \Delta_X \]
holds on $X \times X$, see for example \cite{fay}, p. 31. Denote
by $p_1 : X \times X \to X$ the projection on the first factor.
Using Proposition \ref{normWeier}, the above equality of divisors
yields a canonical isomorphism of line bundles
\[ \Phi^*(O(\Theta)) \isom p_1^*\left(\Omega_X^{\otimes g(g+1)/2} \otimes \left( \wedge^g \difftials
\otimes_\cc O_X) \right)^\lor \right) \otimes O_{X \times
X}(\Delta_X)^{\otimes g} \] on $X \times X$. We will reprove this
isomorphism in the next proposition, and show that its norm is
constant on $X \times X$. After Corollary \ref{bridge} to this
proposition, the proofs of Theorems \ref{Green} and \ref{delta}
are just a few lines.
\begin{proposition} \label{crux}
On $X \times X$, there exists a canonical
isomorphism of line bundles
\[ \Phi^*(O(\Theta)) \isom p_1^*\left( \Omega_X^{\otimes g(g+1)/2} \otimes \left( \wedge^g \difftials
\otimes_\cc O_X) \right)^\lor \right) \otimes O_{X \times
X}(\Delta_X)^{\otimes g} \, . \] The norm of this isomorphism is
everywhere equal to $\exp(\delta(X)/8)$.
\end{proposition}
\begin{proof} We are done if we can prove that
\[ \exp(\delta(X)/8) \cdot \| \vartheta \|(gP-Q) = \|
\tilde{\omega} \|_{\mathrm{Ar}}(P) \cdot G(P,Q)^g \] for all $P,Q
\in X$, where $\{ \omega_1,\ldots,\omega_g \}$ is an orthonormal
basis of $\difftials$. But this follows from the formula in
Theorem \ref{Faltings}, by a computation which is performed in
\cite{jo}, p.~233. Let $P$ be a point on $X$, and choose a local
coordinate $z$ about $P$. By definition of the canonical metric on
$\Omega_X^1$ we have then that $\lim_{Q \to P} |z(Q)-z(P)|/G(Q,P)
= \|dz\|_{\mathrm{Ar}}(P)$. Letting $P_1,\ldots,P_g$ approach $P$
in Theorem \ref{Faltings} we get
 \[
\begin{aligned}
\lim_{P_l \to P} \frac{ \| \det \omega_k(P_l) \|_{\mathrm{Ar}} }{
\prod_{k < l} G(P_k,P_l) }
 &=  \lim_{P_l \to P} \left\{
  \frac{ \| \det \omega_k(P_l) \|_{\mathrm{Ar}} }{ \prod_{k<l} |z(P_k)-z(P_l)|} \cdot
  \frac{\prod_{k<l} |z(P_k)-z(P_l)|}{\prod_{k < l} G(P_k,P_l)} \right\} \\
  &=  \left\{ \lim_{P_l \to P}
      \frac{ | \det \omega_k(P_l) | }{ \prod_{k<l} |z(P_k)-z(P_l)|}  \right\}
      \cdot \| dz \|_{\mathrm{Ar}}^{g + g(g-1)/2}(P) \\
  &=  | W_z(\omega)(P) | \cdot \| dz \|_{\mathrm{Ar}}^{g(g+1)/2}(P) \\
  &=  \| \tilde{\omega} \|_{\mathrm{Ar}}(P) \, .
\end{aligned} \]
The required formula is therefore just a limiting case of Theorem
\ref{Faltings} where all $P_k$ approach $P$.
\end{proof}
\begin{corollary} \label{RandS}
The formula $S(X) = R(X) \cdot
\exp(\delta(X)/8)$ holds.
\end{corollary}
\begin{proof} This follows easily by taking logarithms in the
formula \[ \exp(\delta(X)/8) \cdot \| \vartheta \|(gP-Q) = \|
\tilde{\omega} \|_{\mathrm{Ar}}(P) \cdot G(P,Q)^g \] and
integrating against $\mu(P)$. Here we use again property (iii) of the
Arakelov-Green function and the formula $\log R(X) = - \int_X \log
\|\tilde{\omega}\|_{\mathrm{Ar}}(P) \cdot \mu(P)$, which was noted
above.
\end{proof}
\begin{corollary} \label{bridge}
(1) Let $Q \in X$. Then we have a canonical
isomorphism
\[ \phi_Q^*(O(\Theta)) \isom O_X(\ww + g \cdot Q) \]
of constant norm $S(X)$ on $X$. (2) We have a canonical
isomorphism
\[ \left( \Phi^*(O(\Theta))|_{\Delta_X} \right) \otimes_{O_X}
\Omega_X^{\otimes g} \isom O_X(\ww) \] of constant norm $S(X)$ on
$X$.
\end{corollary}
\begin{proof} We obtain the isomorphism in (1) by restricting the
isomorphism from Proposition \ref{crux} to a slice $X \times \{ Q
\}$, and using Proposition \ref{normWeier}. Its norm is then equal
to $R(X) \cdot \exp(\delta(X)/8)$, which is $S(X)$ by Corollary
\ref{RandS}. For the isomorphism in (2) we restrict the
isomorphism from Proposition \ref{crux} to the diagonal, and apply
the canonical adjunction isomorphism $O_{X \times
X}(-\Delta_X)|_{\Delta_X} \isom \Omega_X^1$. Again we get norm
equal to $R(X) \cdot \exp(\delta(X)/8)$, since the adjunction
isomorphism is an isometry.
\end{proof}
Note that Corollary \ref{bridge} gives an alternative interpretation to
the invariant $S(X)$.
\begin{proof} [Proof of Theorem \ref{Green}] By taking norms of
canonical sections on left and right in the isomorphism from
Corollary \ref{bridge} (1) we obtain
\[ G(P,Q)^g \cdot
\prod_{W \in \ww} G(P,W) = S(X) \cdot \| \vartheta \| (gP-Q) \]
for any $P,Q \in X$. Now take the (weighted) product over $Q \in
\ww$. This gives \[ \prod_{W \in \ww} G(P,W)^{g^3} = S(X)^{g^3-g} \cdot
\prod_{W \in \ww} \| \vartheta \|(gP-W) \, . \] Plugging this in
in the first formula gives
\[ G(P,Q)^g \cdot S(X)^{\frac{g^3-g}{g^3}} \cdot \prod_{W \in \ww} \| \vartheta
\|(gP-W)^{1/g^3} = S(X) \cdot \| \vartheta \|(gP-Q) \, ,
\] from which the theorem follows.
\end{proof}
\begin{proof} [Proof of Theorem \ref{delta}] From Corollary
\ref{bridge} (1) we obtain, again by taking the (weighted) product
over $Q \in \ww$, a canonical isomorphism
\[ \left( \bigotimes_{W \in \ww} \phi_W^*O(\Theta) \right) \isom
O_X(g^3 \cdot \ww) \] of norm $S(X)^{g^3-g}$. It follows that we
have a canonical isomorphism
 \[ \left( \bigotimes_{W \in \ww} \phi_W^*O(\Theta)
 \right)^{\otimes (g-1)/g^3} \isom
O_X((g-1) \cdot \ww) \] of norm $S(X)^{(g-1)(g^3-g)/g^3}$. From
Corollary \ref{bridge} (2) we obtain a canonical isomorphism
\[ \left( \left( \Phi^*(O(\Theta))|_{\Delta_X} \right)
\otimes_{O_X} \Omega_X^{\otimes g} \right)^{\otimes -(g+1)} \isom
O_X(-(g+1)\ww)
\] of norm $S(X)^{-(g+1)}$. Finally from Proposition \ref{normWeier}
and Corollary \ref{RandS} we have a canonical isomorphism
\[  \left( \Omega_X^{\otimes g(g+1)/2} \otimes_{O_X} \left( \wedge^g \difftials
\otimes_\cc O_X) \right)^\lor \right)^{\otimes 2} \isom O_X(2\ww)
\] of norm $S(X)^2 \exp(-\delta(X)/4)$. It follows that indeed the
line bundle $L_X$ is canonically trivial, and that its canonical
trivialising section has norm
\[ S(X)^{-(g-1)(g^3-g)/g^3} \cdot S(X)^{g+1} \cdot S(X)^{-2} \cdot
\exp(\delta(X)/4) = S(X)^{(g-1)/g^2} \cdot \exp(\delta(X)/4)
\, . \] By definition this is $T(X)$, so the
theorem follows.
\end{proof}
It remains to make clear that the invariant $T(X)$ admits an elementary
description in terms of classical functions.
\begin{proposition} \label{classicalTX}
Let $P \in X$ not a Weierstrass point and
let $z$ be a local coordinate about $P$. Define $\|F_z\|(P)$ as
\[ \|F_z\|(P) := \lim_{Q \to P} \frac{ \|\vartheta\|(gP-Q)
}{|z(P)-z(Q)|^g} \, . \] Let $\{\omega_1,\ldots,\omega_g \}$ be an
orthonormal basis of $\difftials$. Then the formula
\[ T(X) = \| F_z \|(P)^{-(g+1)} \cdot \prod_{W
\in \ww} \| \vartheta \|(gP-W)^{(g-1)/g^3} \cdot
|W_z(\omega)(P)|^2 \] holds.
\end{proposition}
\begin{proof} Let $F$ be the canonical section of $\left( \Phi^*(O(\Theta))|_{\Delta_X} \right) \otimes
\Omega_X^{\otimes g}$ given by the canonical isomorphism in
Corollary \ref{bridge} (2). For its norm we have $\| F \|= \|F_z
\| \cdot \|dz\|_{\mathrm{Ar}}^g$ in the local coordinate $z$. The
canonical section of $ \bigotimes_{W \in \ww} \phi_W^*O(\Theta) $
has norm $\prod_{W \in \ww} \| \vartheta \|(gP-W)$ at $P$.
Finally, the canonical section of $\Omega_X^{\otimes g(g+1)/2}
\otimes_{O_X} \left( \wedge^g \difftials \otimes_\cc O_X)
\right)^\lor$ has norm $\| \tilde{\omega} \|_{\mathrm{Ar}} =
|W_z(\omega)| \cdot \|dz\|_{\mathrm{Ar}}^{g(g+1)/2}$. The
proposition follows then from the definition of $T(X)$.
\end{proof}
In \cite{gu}, Gu\`ardia introduced a function $\|J\|$ on
$\mathrm{Sym}^g X$ which involves the first order partial
derivatives of the theta function. We claim that it can be used to
give a formula for $T(X)$ which is especially well-suited for
concrete calculations. Let $\tau \in \hh_g$, the Siegel upper half
space of complex symmetric $g \times g$-matrices with positive
definite imaginary part, be a period matrix associated to $X$.
Consider then the analytic jacobian $\mathrm{Jac}(X) := \cc^g
/\zz^g + \tau  \zz^g$. Then for $w_1,\ldots,w_g \in \cc^g$ we
put
\[ \begin{array}{rcl}
J(w_1,\ldots,w_g) & := & \det \left( \frac{\partial \vartheta}{\partial
z_k}(w_l) \right) \, , \\
\|J\|(w_1,\ldots,w_g) & := & (\det \imtau)^{\frac{g+2}{4}} \cdot \exp( - \pi
\sum_{k=1}^g {}^t y_k (\imtau)^{-1} y_k ) \cdot |J(w_1,\ldots,w_g)| \,
. \end{array} \]
Here $y_k = \mathrm{Im} w_k$ for $k=1,\ldots,g$. The latter definition depends only on the classes in $\mathrm{Jac}(X)$
of the vectors $w_k$.
For a set of $g$ points $P_1,\ldots,P_g$ on $X$ we let, under the
usual correspondence $\mathrm{Pic}_{g-1}(X) \leftrightarrow
\mathrm{Jac}(X)$, the divisor $\sum_{ l=1 \atop l \neq k}^g P_l$
correspond to the class $[w_k] \in \mathrm{Jac}(X)$ of a vector $w_k
\in \cc^g$. We then define $\|J\|(P_1,\ldots,P_g) :=
\|J\|(w_1,\ldots,w_g)$; one may check that this does not depend on the
choice of the period matrix $\tau$ at the beginning.
The following theorem is
Corollary 2.6 in \cite{gu}.
\begin{theorem} \label{guardia} Let $P_1,\ldots,P_g,Q$
be generic points on $X$. Then the formula
\[ \| \vartheta \| (P_1+ \cdots + P_g - Q)^{g-1} = \exp (\delta(X)/8)
\cdot \|J\|(P_1,\ldots,P_g) \cdot \frac{\prod_{k=1}^{g}
G(P_k,Q)^{g-1} }{ \prod_{k < l} G(P_k,P_l)} \] holds.
\end{theorem}
\begin{proposition} \label{Twithfirstorder}
Let $P_1,\ldots,P_g,Q$
be generic points on $X$. Then the formula
\[ \begin{aligned}
T(X) &= \left( \frac{ \| \vartheta \|(P_1 +\cdots+ P_g-Q)}{
\prod_{k=1}^g \| \vartheta \|(g P_k - Q)^{1/g}} \right)^{2g-2}
\cdot
\\ &\quad \quad \cdot \left( \frac{ \prod_{k \neq l} \| \vartheta \|(g
P_k - P_l)^{1/g}}{ \|J\|(P_1,\ldots,P_g)^2 } \right) \cdot
\prod_{W \in \ww} \prod_{k=1}^g \| \vartheta \|(g P_k -
W)^{(g-1)/g^4}  \end{aligned} \] holds.
\end{proposition}
\begin{proof} The formula follows from Theorem \ref{guardia},
using Theorem \ref{Green} to eliminate the
occurring values of the Arakelov-Green function $G$,
and using Theorem \ref{delta} to eliminate the factor
$\exp(\delta(X)/8)$. The factors involving $S(X)$ that are
introduced in this way cancel out.
\end{proof}
For example, if $g=1$ and $X$ is given as $X =\cc/\zz+\tau \zz$
with $\tau$ in the complex upper half plane, we obtain \[
T(X) = (\imtau)^{-3/2} \exp(\pi \imtau/2) \cdot | \frac{\partial
\vartheta}{\partial z} ( \frac{1+\tau}{2};\tau )|^{-2}
\, . \] By Jacobi's derivative formula we have then
\[ T(X) = (2\pi)^{-2} \cdot ((\imtau)^6 |\Delta(\tau)|)^{-1/4} \]
where $\Delta(\tau)$ is the discriminant modular form
$\Delta(q):=\eta(q)^{24} = q \prod_{n=1}^\infty (1-q^n)^{24}$. It
follows that Faltings' delta-invariant is given by
\[ \delta(X) = -\log((\imtau)^6 |\Delta(\tau)|)-8 \log(2\pi) \]
which is well-known, see \cite{fa}, p. 417.

The formula for $T(X)$ for an elliptic curve $X$ can be
generalised to hyperelliptic Riemann surfaces of arbitrary genus.
In \cite{jong3} the following result is proven. For any integer $g
\geq 2$, let $\varphi_g$ be the discriminant modular form on
$\hh_g$ as defined in \cite{lock}, Section 3. This is a modular
form on $\Gamma_g(2) := \{ \gamma \in \mathrm{Sp}(2g,\zz) \, : \,
\gamma \equiv I_{2g} \, \bmod 2 \}$ of weight $4r$, where
$r:={2g+1 \choose g+1}$.
\begin{theorem} \label{hyperelliptic}
Let $X$ be a hyperelliptic Riemann surface of genus $g \geq 2$.
Choose an ordering of the Weierstrass points on $X$ and a
canonical symplectic basis of the homology of $X$ given by this
ordering (cf. \cite{mu}, Chapter IIIa, \S 5). Let $\tau \in \hh_g$
be the period matrix of $X$ associated to this canonical basis and
put $\Delta_g(\tau) := 2^{-(4g+4)n} \cdot \varphi_g(\tau)$ where
$n := {2g \choose g+1}$. Then the formula
\[ T(X) = (2\pi)^{-2g} \cdot ((\imtau)^{2r}
|\Delta_g(\tau)|)^{-\frac{3g-1}{8ng}} \] holds.
\end{theorem}
The proof of Theorem \ref{hyperelliptic} is quite complicated, and
unfortunately we do not know how to generalise the proof to
arbitrary Riemann surfaces of genus $g$. We leave it as an open
question whether in general the invariant $T(X)$ can be naturally
expressed in terms of Siegel modular forms on $\hh_g$.

\section{Applications to intersection theory} \label{application}

In this section we use Proposition \ref{normWeier} to give a
formula for the relative dualising sheaf on a semi-stable
arithmetic surface (Proposition \ref{omega}). As consequences we
derive a lower bound for the self-intersection of the relative
dualising sheaf (Proposition \ref{omegasq}) and a formula for the
self-intersection of a point (Proposition \ref{selfinterspoint}).

Let $p: \xx \to B$ be a semi-stable arithmetic surface over the
spectrum $B$ of the ring of integers in a number field $K$. We
assume that the generic fiber $\xx_K$ is a geometrically
connected, smooth proper curve of genus $g>0$. Denote by $\ww$ the
Zariski closure in $\xx$ of the divisor of Weierstrass points on
$\xx_K$, and denote by $\omega_{\xx/B}$ the relative dualising
sheaf of $p$. 

The next lemma is an analogue of Lemma 3.3 in \cite{ar2}.
\begin{lemma} \label{withoutnorms1}
There exists an effective vertical divisor $V$ on $\xx$ such that
we have a canonical isomorphism
\[  \omega_{\xx/B}^{\otimes g(g+1)/2 } \otimes_{O_\xx} \left(p^*(\det p_*
\dual) \right)^\lor {\buildrel \sim \over \longrightarrow} O_\xx(V+\ww)
 \] of line bundles on $\xx$.
\end{lemma}
\begin{proof} We have on $\xx$ a
canonical sheaf morphism $ p^*(\det p_* \dual) \longrightarrow
\omega_{\xx/B}^{\otimes g(g+1)/2 } $ given locally by \[ \xi_1
\wedge \ldots \wedge \xi_g \mapsto \frac{\xi_1 \wedge \ldots
\wedge \xi_g}{\psi_1 \wedge \ldots \wedge \psi_g } \cdot
\tilde{\psi}  \] for a $K$-basis $\{ \psi_1,\ldots,\psi_g \}$ of
the differentials on the generic fiber of $\xx$. Multiplying by $
(p^*(\det p_* \dual))^\lor $ we obtain a morphism
\[ O_\xx \longrightarrow \omega_{\xx/B}^{\otimes g(g+1)/2 }
\otimes_{O_\xx} \left( p^*(\det p_* \dual) \right)^\lor \, . \]
The image of 1 is a section whose divisor is an effective divisor
$V+\ww$ where $V$ is vertical. This gives the required
isomorphism.
\end{proof}
We will now turn to the Arakelov intersection theory on $\xx$. Our
references are, once more, \cite{ar} and \cite{fa}. For a complex
embedding $\sigma : K \hookrightarrow \cc$ we denote by $F_\sigma$
the ``fiber at infinity'' associated to $\sigma$. The
corresponding Riemann surface of genus $g$
is denoted by $X_\sigma$.
\begin{proposition} \label{omega} Let $V$ be the effective vertical divisor from
Lemma \ref{withoutnorms1}. Then we have
\[ \frac{1}{2} g (g+1) \dual = V + \ww + \sum_{\sigma : K \hookrightarrow \cc} \log R(X_\sigma)
\cdot F_\sigma + p^*(\det p_* \omega_{\xx/B})   \] as Arakelov
divisors on $\xx$. Here the sum runs over the embeddings
of $K$ in $\cc$.
\end{proposition}
\begin{proof} Consider the canonical isomorphism from Lemma
\ref{withoutnorms1}. The restriction of this isomorphism to
$X_\sigma$ is the isomorphism of Proposition \ref{normWeier}. In
particular it has norm $R(X_\sigma)$. The proposition follows.
\end{proof}
We shall deduce two consequences from this proposition. We assume
for the moment that $g \geq 2$. We define $R_b$ for a closed point
$b \in B$ by the equation $ (2g-2) \cdot \log R_b = (V_b,\dual) $,
where the intersection is taken in the sense of Arakelov.
The assumption that $p: \xx \to B$ is semi-stable implies that the quantity
$\log R_b $ is always non-negative.
\begin{proposition} \label{omegasq} Assume that $g \geq 2$. Then
the lower bound
\[ (\dual,\dual) \geq \frac{ 8(g-1) }{ (2g-1)(g+1) } \cdot \left(
\sum_b \log R_b + \sum_{\sigma : K \hookrightarrow \cc} \log
R(X_\sigma) + \height \right)
\] holds. Here the first sum runs over the closed points $b \in
B$, and the second sum runs over the embeddings of $K$ in
$\cc$.
\end{proposition}
\begin{proof}
Intersecting the equality from Proposition \ref{omega} with
$\dual$ we obtain
\[ \begin{aligned}
\frac{1}{2}g(g+1) & (\dual,\dual) = \\ &= (\ww,\dual) + (2g-2)
\cdot \left(
  \sum_b \log R_b + \sum_{\sigma : K \hookrightarrow \cc}
  \log R(X_\sigma) + \height \right) \, .
  \end{aligned} \] Now since the generic degree of $\ww$ is
$g^3-g$ we obtain by Theorem 5 of \cite{fa} the lower bound
\[ (\ww,\dual) \geq
\frac{g^3-g}{2g(2g-2)} (\dual,\dual) \, . \] Using this in the
first equality
gives the result.
\end{proof}
One should compare the above lower bound for $(\dual,\dual)$ with
the lower bounds for $(\dual,\dual)$ given in \cite{bu}, Section
3.3. The contributions at infinity $\log R(X_\sigma)$ have
properties similar to the terms $A_{k,\sigma}$ occurring in
\cite{bu}. In particular, the right-hand side of the inequality
in Proposition \ref{omegasq} may be negative.

We refer to the author's thesis for a proof of the
following result.
\begin{proposition} \label{asymptR}
Let $X_t$ be a holomorphic family of compact and connected Riemann
surfaces of genus $g \geq 2$ over the punctured disk, degenerating
to the union of two Riemann surfaces of positive genera $g_1$,
$g_2$ with two points identified. Suppose that neither of these
two points was a Weierstrass point on each of the two separate
Riemann surfaces. Then the formula
\[ \log R(X_t) = -\frac{g_1g_2}{2g} \log |t| + O(1) \quad \mathrm{as}
\quad t \to 0   \] holds.
\end{proposition}
In particular, the value $\log R(X_t)$ goes to plus infinity under
the conditions described in the theorem. It would be interesting
to have a more precise, quantitative version of Proposition
\ref{asymptR}.

Our second result is a formula for the self-intersection of a point.
In the proof of the next lemma, we make use of the Deligne bracket
(see \cite{de}). This is a rule that assigns to a pair $L,M$ of
line bundles on $\xx$ a line bundle $\langle L,M \rangle$ on $B$
such that the following properties hold: (i) we have canonical
isomorphisms $ \langle L_1 \otimes L_2,M \rangle \isom \langle
L_1,M \rangle \otimes \langle L_2,M \rangle$, $\langle L,M_1
\otimes M_2 \rangle \isom \langle L,M_1 \rangle \otimes \langle
L,M_2 \rangle$ and $\langle L,M \rangle \isom \langle M,L
\rangle$; (ii) for a section $P: B \to \xx$ we have a canonical
isomorphism $\langle O_\xx(P),M \rangle \isom P^*M$; (iii)
(adjunction formula) for a section $P:B \to \xx$ we have a
canonical isomorphism $\langle P,\omega_{\xx/B} \rangle \isom
\langle P,P \rangle^{\otimes -1}$; (iv) (Riemann-Roch) for a line
bundle $L$ on $\xx$ we have a canonical isomorphism $(\det
Rp_*L)^{\otimes 2} \isom \langle L,L \otimes \omega_{\xx/B}^{ -1}
\rangle \otimes (\det p_* \omega)^{\otimes 2}$ relating the
Deligne bracket to the determinant of cohomology. 

Assume that $g \geq 1$ again.
\begin{lemma} \label{withoutnorms2}
Let $P$ be a section of $p$, not a Weierstrass point
on the generic fiber. Then we have a canonical isomorphism
\[ P^*( O_\xx(V + \ww))^{\otimes 2} \isom
\left( \det R  p_* O_\xx(g  P) \right)^{\otimes 2} \] of
line bundles on $B$.
\end{lemma}
\begin{proof}  Applying Riemann-Roch to the line bundle
$O_\xx(g  P)$ we obtain a canonical isomorphism \[ \left( \det Rp_*
O_\xx(g  P) \right)^{\otimes 2} \isom \langle O_\xx(g  P),
O_\xx( g P) \otimes \omega_{\xx/B}^{-1} \rangle \otimes
(\det p_* \omega_{\xx/B})^{\otimes 2}  \] of line bundles on $B$.
The line bundle at the right hand side is, by the adjunction
formula, canonically isomorphic to the line bundle $\langle P,P
\rangle^{\otimes g(g+1) } \otimes (\det p_*
\omega_{\xx/B})^{\otimes 2}$. On the other hand, pulling back the
isomorphism from Lemma \ref{withoutnorms1} along $P$ and using once more the adjunction
formula gives a canonical isomorphism \[ \langle P,P
\rangle^{\otimes -g(g+1)/2} \isom \langle V+\ww, P \rangle \otimes
\det p_* \omega_{\xx/B} \, . \] The lemma follows by a
combination of these observations.
\end{proof}
\begin{proposition} \label{selfinterspoint}
Let $P$ be a section of $p$, not a Weierstrass
point on the generic fiber. Then $-\frac{1}{2}g(g+1)  (P,P)$ is
given by the expression
\[  - \sum_{\sigma : K
\hookrightarrow \cc} \log G(P_\sigma,\ww_\sigma) + \log \# R^1 p_*
O_\xx(g \cdot P) + \sum_{\sigma : K \hookrightarrow \cc} \log
R(X_\sigma) + \height \, ,
\] where $\sigma$ runs through the complex embeddings of $K$.
\end{proposition}
\begin{proof} Intersecting the equality from Proposition
\ref{omega} with $P$, and using the adjunction formula
$(\omega,P)=-(P,P)$, we obtain the equality
\[ -\frac{1}{2}g(g+1)  (P,P) = (V+\ww,P) + \sum_{\sigma : K
\hookrightarrow \cc} \log R(X_\sigma) + \height \, . \] It remains
therefore to see that $(V+\ww,P)_{\mathrm{fin}}=\log \# R^1 p_*
O_\xx(g \cdot P)$. For this we consider the isomorphism in Lemma
\ref{withoutnorms2}. Note that $p_* O_\xx(g \cdot P)$ is
canonically trivialised by the function 1.
This gives a canonical section at the right hand side with norm the
square of $\# R^1 p_*O_\xx(g \cdot P)$. Under the isomorphism, it
is identified with the canonical section on the left-hand side,
which has norm the square of $\exp((V+\ww,P)_\mathrm{fin})$.
The required equality follows.
\end{proof}
We see that minus the self-intersection of a point $P$ is large if
$P$ is close to a Weierstrass point, either in the $p$-adic or in
the complex topology.

\section{A numerical example} \label{example}

In this final section we wish to illustrate the practical
significance of our Theorems \ref{Green} and \ref{delta} by
exhibiting a concrete example dealing with a hyperelliptic curve
of genus 3. The propositions below can be proved by methods
similar to those in \cite{bmmb}, Section 3.

Let $K$ be a number field, and $A$ its ring of
integers. Let $F \in A[x]$ be monic of degree 5 with $F(0)$ and
$F(1)$ a unit in $A$. Put $R(x) := x(x-1) + 4 F(x)$. Suppose that the
following holds for $R$: the discriminant $\Delta$ of $R$ is non-zero; for
every prime $\wp$ of residue characteristic $\mathrm{char}(\wp) \neq 2$ of $A$ we have
$v_\wp(\Delta)=0$ or 1; if $\mathrm{char}(\wp) \neq 2$ and
$v_\wp(\Delta)=1$, then $R (\bmod \wp)$ has a unique multiple root,
and its multiplicity is 2.
\begin{proposition} \label{geometry}
The equation
\[ C_F \, : \, y^2 = x(x-1)R(x) \] defines a hyperelliptic curve
of genus 3 over $K$. It extends to a semi-stable arithmetic
surface $p: \xx \to B=\mathrm{Spec}(A)$. We have that $\xx$ has
bad reduction at $\wp$ if and only if $\mathrm{char}(\wp) \neq 2$
and $v_\wp(\Delta)=1$. In this case, the bad fiber is an irreducible
curve with a single double point. The differentials
$dx/y,xdx/y,x^2dx/y$ form a basis of the $O_B$-module $p_*
\omega_{\xx/B}$. The points $W_0,W_1$ on $C_F$ given by $x=0$ and
$x=1$ extend to disjoint $\sigma$-invariant sections of $p$.
\end{proposition}
As for the Arakelov invariants of $C_F$, we have the following
result.
\begin{proposition} \label{concrete} At a complex embedding
$\sigma : K \hookrightarrow \cc$, let $\Omega_\sigma =
(\Omega_{1\sigma} | \Omega_{2 \sigma})$ be a period matrix for the
Riemann surface corresponding to $C_F \otimes_{\sigma,K} \cc$,
formed on the basis $dx/y,xdx/y,x^2dx/y$. Further, let
$\tau_\sigma = \Omega_{1\sigma}^{-1} \Omega_{2\sigma}$. Then
\[ \Ardeg \det p_* \omega_{\xx/B} = -\frac{1}{2} \sum_\sigma \log
\left( |\det \Omega_{1 \sigma}|^2 (\det \imtau_\sigma) \right) \,
,
\] where
the sum runs over the complex embeddings of $K$. Further, the
formula \[ (\omega_{\xx/B}, \omega_{\xx/B}) = 24 \sum_\sigma
\log G_\sigma(W_0,W_1) \] holds.
\end{proposition}

For our example, we choose the polynomial
$F(x)=x^5+6x^4+4x^3-6x^2-5x-1$ defined over $\qq$. Then the
corresponding $R(x) = x(x-1)+4F(x)$ satisfies the conditions
described above. The corresponding hyperelliptic curve (which we
will call $X$ from now on) of genus 3 has bad reduction at the
primes $p=37, p=701$ and $p=14717$. An equation is given by
\[ X \, : \, y^2 = x(x-1)(4x^5+24x^4+16x^3-23x^2-21x-4) \, . \]
We choose an ordering of the Weierstrass points of $X$. As in
\cite{mu}, Chapter III, \S 5 we construct from this a canonical
symplectic basis of the homology of (the Riemann surface
corresponding to) $X$. Using Mathematica, we compute the
periods of the differentials $dx/y, xdx/y, x^2
dx/y$. This leads to an explicit value of
$\Omega=(\Omega_1|\Omega_2)$ and
the numerical approximation
\[ \Ardeg \det p_* \omega_{\xx/B} = -1.280295247656532068... \]
Using the Riemann vector given by \cite{mu}, p. 3.82 we can make the
identification $\mathrm{Pic}_2(X) \leftrightarrow \cc^3/\zz^3 +
\tau \zz^3$ explicit. With Theorem \ref{hyperelliptic} we find then the
following numerical approximation to $T(X)$:
\[ \log T(X) = -4.44361200473681284... \]
The values of the theta function that are needed for this computation
are approximated by the defining summation
formula, which consists of rapidly decreasing
exponential terms. An elementary \emph{a priori} calculation shows
how much terms we need to compute in order to approximate a value of
the theta function with a prescribed accuracy.

It remains then to calculate the invariant $\log S(X)$. Recall the
definition
\[ \log S(X) := - \int_X \log \|\vartheta\|(3P-Q) \cdot \mu(P) \, . \]
Note that the integrand diverges at infinity, so we would rather
want to make use of the formula
\[ \log S(X) = -9 \int_X \log \|\vartheta\|(3P-Q) \cdot \mu(Q)
+ \frac{1}{3} \cdot \sum_{W \in \ww} \log \|\vartheta\|(3P-W) \, , \]
valid for any $P \in X$ which is not a Weierstrass point. This
formula can be easily derived from Theorem
\ref{Green} by taking logarithms and integrating
against $\mu(Q)$. The integrand has now only a (logarithmic)
singularity at $Q=P$.
Write $x=u + iv$ with $u,v$ real. We want to express $\mu(Q)$ in
terms of the coordinates $u,v$. This is done by the following lemma.
\begin{lemma} Let $h$ be the $3 \times 3$-matrix given by
\[ h = \left( \overline{\Omega}_1 (\imtau) {}^t \Omega_1 \right)^{-1}
\, . \] Then we can write
\[ \mu = \left( h_{11} + 2h_{12}u + 2h_{13}(u^2-v^2) + h_{22}(u^2+v^2)
+ 2h_{23} u(u^2+v^2) + h_{33} (u^2+v^2)^2 \right) \cdot \frac{du dv}{
3|f| } \] in the coordinates $u,v$.
\end{lemma}
\begin{proof} Let $\omega_k = x^{k-1}dx/y$ for $k=1,2,3$. By
Riemann's bilinear relations, the
fundamental (1,1)-form $\mu$ is given by $\mu = \frac{i}{6} \sum_{k,l=1}^3
h_{kl} \cdot \omega_k \wedge \overline{\omega_l} $. Expanding this
expression gives the result, where we note that the matrix $h$ is real
symmetric, since our defining equation for $X$ is defined over the
reals.
\end{proof}
We can now effectuate the required integral, choosing an arbitrary point
$P$ and taking care of the singularity of the integrand at
this point $P$. We find the approximation
\[ \log S(X) = 17.57... \]
In order to check this result, we have taken several choices for $P$.
By Theorem \ref{delta} we have
\[ \delta(X) = -33.40... \] and using Theorem \ref{Green} we can
approximate, by taking $Q=W_1$ and letting $P$ approach $W_0$,
\[ G(W_0,W_1) = 2.33... \]
By Proposition \ref{concrete} we finally find
\[ (\omega_{\xx/B}, \omega_{\xx/B}) = 20.32... \]
The running times of the computations were negligible, except for the
computation of the integral involved in
$\log S(X)$, which took about 7 hours on the
author's laptop.

\subsection*{Acknowledgments} The author wishes to thank his thesis
advisor Gerard van der Geer for his encouragement and helpful
remarks.

\end{document}